\documentstyle[thmsa,11pt]{article}

\renewcommand{\abstract}{ABSTRACT.}
\newtheorem{Theorem}{\bf Theorem}[section]
\newtheorem{Proposition}[Theorem]{ \bf Proposition}
\newtheorem{Corollary}[Theorem]{ \bf Corollary}
\newtheorem{Definition}[Theorem]{\it  Definition}
\newtheorem{Lemma}[Theorem]{\bf Lemma}
\newtheorem{Remark}[Theorem]{\it Remark}
\input{tcilatex}
\input tcilatex
\begin{document}

\title{CROSSED\ PRODUCTS\ OF\ LOCALLY\ \ $C^{*}$-ALGEBRAS}
\author{MARIA\ JOI\c{T}A}
\maketitle

\begin{abstract}
The crossed products of locally $C^{*}$-algebras are defined and a Takai
duality theorem for inverse limit actions of a locally compact group on a
locally $C^{*}$-algebra is proved.

\end{abstract}

\section{ Introduction}

Locally $C^{*}$-algebras are generalizations of $C^{*}$-algebras. Instead of
being given by a single $C^{*}$-norm, the topology on a locally $C^{*}$%
-algebra is defined by a directed family of $C^{*}$-semi-norms. In [{\bf 9}%
], Phillips defines the notion of action of a locally compact group $G$ on a
locally $C^{*}$-algebra $A$ whose topology is determined by a countable
family of $C^{*}$-semi-norms, and also defines the crossed product of $A$ by
an inverse limit action $\alpha =\lim\limits_{\stackunder{n}{\leftarrow }%
}\alpha ^{(n)}$ as being the inverse limit of crossed products of $A_n$ by $%
\alpha ^{(n)}$. In this paper, by analogy with the case of $C^{*}$-algebras,
we define the concept of crossed product, respectively reduced crossed
product of locally $C^{*}$-algebras.

The Takai duality theorem says that if $\alpha $ is a continuous action of
an abelian locally compact group $G$ on a $C^{*}$-algebra $A$, then we can
recover the system $(G,A,\alpha )$ up to stable isomorphism from the double
dual system in which $G=\widehat{\widehat{G}}$ acts on the crossed product $%
\left( A\times _{\alpha _{}}G\right) \times _{\widehat{\alpha }}\widehat{G}$
by the dual action of the dual group. In [{\bf 3}], Imai and Takai prove a
duality theorem for $C^{*}$-crossed products by a locally compact group that
generalizes the Takai duality theorem [{\bf 12}]. For a given $C^{*}$%
-dynamical system $(G,A,\alpha ),$ they construct a ''dual '' $C^{*}$%
-crossed product of the reduced crossed product $A\times _{\alpha ,r}G$ by
an isomorphism $\beta $ from $A\times _{\alpha ,r}G$ into $L(H),$ the $C^{*}$%
-algebra of all bounded linear operators on some Hilbert space $H,$ and show
that this is isomorphic to the tensor product $A\otimes {\cal K}(L^2(G))$ of 
$A$ and ${\cal K}(L^2(G)),$ the $C^{*}$-algebra of all compact operators on $%
L^2(G)$. If $G$ is commutative, the ''dual '' $C^{*}$-crossed product
constructed by Imai and Takai is isomorphic to the double crossed product $%
\left( A\times _{\alpha _{}}G\right) \times _{\widehat{\alpha }}\widehat{G}$%
. Katayama [{\bf 6}] shows that a non-degenerate coaction $\beta $ of a
locally compact group on a $C^{*}$-algebra $A$ induces an action $\widehat{%
\beta }$ of $G$ on the crossed product $A\times _{\beta _{}}G$ and proves
that the $C^{*}$-algebras $\left( A\times _{\beta _{}}G\right) \times _{%
\widehat{\beta },r}G$ and $A\otimes {\cal K}(L^2(G))$ are isomorphic. In [%
{\bf 13}], Vallin shows that there is a bijective correspondence between the
set of all actions of a locally compact group $G$ on a $C^{*}$-algebra $A$
and the set of all actions of the commutative Kac $C^{*}$-algebra $C^{*}{\bf %
K}_G^a$ associated with $G$ on $A$. A coaction of $G$ on $A$ is an action of
the symmetric Kac $C^{*}$-algebra $C^{*}{\bf K}_G^s$ associated with $G$. If 
$G$ is commutative, we can identified $C_r^{*}(G)$ with $C_0(\widehat{G})$
via the Fourier transform, whence becomes clear that a coaction of $G$ is
the same thing as an action of $\widehat{G}$. Thus we can regard the
coactions of a locally compact group $G$ as ''actions of the dual group even
there isn't any dual group''. Also, Vallin shows that an action $\alpha $ (
coaction $\beta $ ) of $G$ on $A$ induces a coaction $\widehat{\alpha }$ (
action \ \ $\widehat{\beta }$ ) of $G$ on the crossed product $A\times
_{\alpha ,r}G$ (respectively $A\times _{\beta _{}}G$) and proves a version
of the Takai duality theorem showing that the double crossed product $\left(
A\times _{\alpha ,r}G\right) \times _{\widehat{\alpha }}G$ is isomorphic to $%
A\otimes {\cal K}(L^2(G)).$ We propose to prove a version of the Takai
duality theorem for crossed products of locally $C^{*}$-algebras.

The paper is organized follows. In Section 2 we present some basic
definitions and results about locally $C^{*}$-algebras and Kac $C^{*}$%
-algebras. In Section 3 we define the notion of crossed product (reduced
crossed product) of a locally $C^{*}$-algebra $A$ by an inverse limit action 
$\alpha $ of a locally compact group $G$ and prove some basic properties of
these. Section 4 is devoted to actions of a Kac $C^{*}$-algebra on a locally 
$C^{*}$-algebra. We show that there is a bijective correspondence between
the set of all inverse limit actions of a locally compact group $G$ on a
locally $C^{*}$-algebra $A$ and the set of all inverse limit actions of the
commutative Kac $C^{*}$-algebra $C^{*}{\bf K}_G^a$ on $A$, Proposition 4.4.
As a consequence of this result we obtain: for a compact group $G,$ any
action of the Kac $C^{*}$-algebra $C^{*}{\bf K}_G^a$ on $A$ is an inverse
limit of actions of the Kac $C^{*}$-algebras $C^{*}{\bf K}_G^a$ on $A_p$, $%
p\in S(A)$. In Section 5, using the same arguments as in [{\bf 13}], we show
that any inverse limit action $\alpha $ (coaction $\beta $ ) of a locally
compact group $G$ on a locally $C^{*}$-algebra $A$ induces an inverse limit
coaction $\widehat{\alpha }$ ( action $\widehat{\beta }$ ) of $G$ on the
crossed product $A\times _{\alpha ,r}G$ (respectively $A\times _{\beta
_{}}G),$ Proposition 5.5. Finally, we prove that if $\alpha $ is an inverse
limit action of a locally compact group $G$ on a locally $C^{*}$-algebra $A,$
then there is an isomorphism of locally $C^{*}$-algebras from $\left(
A\times _{\alpha ,r}G\right) \times _{\widehat{\alpha }}G$ onto $A\otimes 
{\cal K}(L^2(G))$ and the inverse limit actions $\widehat{\widehat{\alpha }}$
and $\alpha \otimes $ad$\rho $ are equivalent, Theorem 5.6.

\section{Preliminaries}

A locally $C^{*}$-algebra is a complete complex Hausdorff topological $*$
-algebra $A$ whose topology is determined by a family of $C^{*}$-semi-norms,
see [{\bf 1}], [{\bf 2}], [{\bf 4}], [{\bf 9}], [{\bf 10}]. If $S(A)$ is the
set of all continuous $C^{*}$-semi-norms on $A$, then for each $p\in S(A),$ $%
A_p=A/\ker (p)$ is a $C^{*}$-algebra with respect to the norm induced by $p$%
, and $A=\lim\limits_{\stackunder{p\in S(A)}{\leftarrow }}A_p$. The
canonical maps from $A$ onto $A_p,$ $p\in S(A)$ are denoted by $\pi _p$, the
image of $a$ under $\pi _p$ by $a_p,$ and the connecting maps of the inverse
system $\{A_p\}_{p\in S(A)}$ by $\pi _{pq,}$ $p,q\in S(A)$ with $p\geq q$.

A morphism of locally $C^{*}$-algebras is a continuous $*$-morphism $\Phi $
from a locally $C^{*}$-algebra $A$ to a locally $C^{*}$-algebra $B$. An
isomorphism of locally $C^{*}$-algebras is a morphism of locally $C^{*}$%
-algebras which is invertible and its inverse is a morphism of locally $%
C^{*} $-algebras. An $S$ -morphism of locally $C^{*}$-algebras is a morphism 
$\Phi :A\rightarrow M(B)$, where $M(B)$ is the multiplier algebra of $B$,
with the property that for any approximate unit $\{e_i\}_i$ of $A$ the net $%
\{\Phi (e_i)\}_i$ converges to $1$ with respect to the strict topology on $%
M(B).$ If $\Phi :A\rightarrow M(B)$ is an $S$-morphism of locally $C^{*}$%
-algebras, then it extends to a unique morphism $\overline{\Phi }%
:M(A)\rightarrow M(B)$ of locally $C^{*}$-algebras, see [{\bf 5}].

A Kac $C^{*}$-algebra is a quadruple ${\bf K=}(B,d,j,\varphi )$, where $B$
is a $C^{*}$-algebra, $d$ is a comultiplication on $B$, $j$ is a
coinvolution on $B$, and $\varphi $ is a semi-finite, lower semi-continuous,
faithful weight on $B$, see [{\bf 13}].

Let $A$ and $B$ be two locally $C^{*}$-algebras. The injective tensor
product of the locally $C^{*}$-algebras $A$ and $B$ is denoted by $A\otimes
B $, see [{\bf 2}], and the locally $C^{*}$-subalgebra of $M(A\otimes B)$
generated by the elements $x$ in $M(A\otimes B)$ such that $x(1\otimes
B)+(1\otimes B)x\subseteq A\otimes B$ is denoted by $M(A,B)$. If $G$ is a
locally compact group, then $M(A,C_0(G))$ may be identified with the locally 
$C^{*}$-algebra $C_b(G,A)$ of all bounded continuous functions from $G$ to $%
A $.

Let $G$ be a locally compact group. $C^{*}{\bf K}_G^a=\left(
C_0(G),d_G^a,j_G^a,ds\right) $ is the commutative Kac $C^{*}$-algebra
associated with $G$ and $C^{*}{\bf K}_G^s=\left(
C_r^{*}(G),d_G^s,j_G^s,\varphi _G\right) $ is the symmetric Kac $C^{*}$%
-algebra associated with $G$, see [{\bf 13}].

An action of a Kac $C^{*}$-algebra ${\bf K=}(B,d,j,\varphi )$ on a $C^{*}$%
-algebra $A$ is an injective $S$ -morphism $\alpha $ from $A$ to $M(A,B)$
such that $(\overline{\alpha \otimes \text{id}})\circ \alpha =(\overline{%
\text{id}_A\otimes \sigma _B\circ d})\circ \alpha $, see [{\bf 13}].

\smallskip

\section{Crossed products}

Let $A$ be a locally $C^{*}$-algebra and let $G$ be a locally compact group.

\smallskip

\begin{Definition}
{\em An action of }$G${\em \ on }$A${\em \ is a morphism }$\alpha $ {\em \
from }$G${\em \ to Aut}$\left( A\right) ${\em , the set of all isomorphisms
of locally }$C^{*}${\em -algebras from }$A${\em \ to }$A${\em . The action }$%
\alpha ${\em \ is continuous if the function }$\left( t,a\right) \rightarrow
\alpha _t(a)${\em \ from }$G\times A${\em \ to }$A${\em \ is jointly
continuous.}
\end{Definition}

\smallskip

\begin{Definition}
{\em A locally }$C^{*}${\em -dynamical system is a triple }$\left(
G,A,\alpha \right) ${\em , where }$G${\em \ is a locally compact group, }$A$%
{\em \ is a locally }$C^{*}${\em -algebra and }$\alpha ${\em \ is a
continuous action of }$G${\em \ on }$A.$
\end{Definition}

\smallskip

{\em \smallskip }

\begin{Definition}
{\em We say that }$\left\{ \left( G,A_{\delta _{}},\alpha ^{\left( \delta
\right) }\right) \right\} _{\delta \in \Delta }${\em \ is an inverse system
of }$C^{*}${\em -dynamical systems if }$\left\{ A_{\delta _{}}\right\}
_{\delta \in \Delta }${\em \ is an inverse system of }$C^{*}${\em -algebras
and for each }$t${\em \ in }$G${\em , }$\left\{ \alpha _t^{\left( \delta
\right) }\right\} _{\delta \in \Delta }${\em \ is an inverse system of }$%
C^{*}${\em -isomorphisms.}

{\em Let }$A=\lim\limits_{\stackunder{\delta \in \Delta }{\leftarrow }%
}A_{\delta _{}}${\em \ and }$\alpha _t=\lim\limits_{\stackunder{\delta \in
\Delta }{\leftarrow }}\alpha _t^{\left( \delta \right) }${\em \ for each }$%
t\in G.${\em \ Then the map }$\alpha :G\rightarrow ${\em Aut}$(A)${\em \
defined by }$\alpha (t)=\alpha _t${\em \ is a continuous action of }$G${\em %
\ on }$A${\em \ and }$(G,A,\alpha )${\em \ is a locally }$C^{*}${\em %
-dynamical system. We say that }$(G,A,\alpha )${\em \ is the inverse limit
of the inverse system of }$C^{*}${\em -dynamical systems }$\left\{ \left(
G,A_{\delta _{}},\alpha ^{\left( \delta \right) }\right) \right\} _{\delta
\in \Delta }${\em .}
\end{Definition}

\smallskip

\begin{Definition}
{\em A continuous action }$\alpha ${\em \ of }$G$ {\em on} $A$ {\em is an
inverse limit action if we can write }$A${\em \ as inverse limit }$%
\lim\limits_{\stackunder{\delta \in \Delta }{\leftarrow }}A_{\delta _{}}$%
{\em \ of }$C^{*}${\em -algebras in such a way that there are actions }$%
\alpha ^{\left( \delta \right) }${\em \ of }$G${\em \ on }$A_{\delta _{}}$%
{\em \ such that }$\alpha _t=\lim\limits_{\stackunder{\delta \in \Delta }{%
\leftarrow }}\alpha _t^{\left( \delta \right) }${\em \ for all }$t${\em \ in 
}$G${\em \ (Definition 5.1,\ [{\bf 9}]).}
\end{Definition}

\smallskip

\begin{Remark}
{\em The action }$\alpha ${\em \ of }$G${\em \ on }$A${\em \ is an inverse
limit action if there is a cofinal subset of }$G${\em -invariant continuous }%
$C^{*}${\em -semi-norms on }$A${\em \ ( a continuous }$C^{*}${\em -semi-norm 
}$p${\em \ on }$A${\em \ is }$G${\em \ -invariant if }$p(\alpha _t(a))=p(a)$%
{\em \ for all }$a${\em \ in }$A${\em \ and for all }$t${\em \ in }$G${\em ).%
}
\end{Remark}

\smallskip

The following lemma is Lemma 5.2 of [{\bf 9}].

\begin{Lemma}
Any continuous action of a compact group $G$\ on a locally $C^{*}$-algebra $%
A $\ is an inverse limit action.
\end{Lemma}

\smallskip

\smallskip {\em \smallskip }Let $(G,A,\alpha )$ be a locally $C^{*}$%
-dynamical system such that $\alpha $ is an inverse limit action. By Remark
3.5, we can suppose that $S(A)$ coincides with the set of all $G$ -invariant
continuous $C^{*}$-semi-norms on $A$.

Let $C_c(G,A)\;$be the vector space of all continuous functions from $G$ to $%
A$ with compact support.

\begin{Lemma}
Let $f\in C_c(G,A)$. Then there is a unique element $\tint\limits_Gf(s)ds$\
in $A$\ such that for any non-degenerate $*$\ -representation $(\varphi
,H_{\varphi ^{}})$\ of $A$\ 
\[
\left\langle \varphi (\tint\limits_Gf(s)ds)\xi ,\eta \right\rangle
=\tint\limits_G\left\langle \varphi (f(s))\xi ,\eta \right\rangle ds 
\]
for all $\xi ,\eta $\ in $H_{\varphi ^{}}$. Moreover, we have:

\begin{quote}
(1) $p(\tint\limits_Gf(s)ds)\leq M\sup \{p(f(s));$\ $s\in $supp$\left(
f\right) \}$\ for some positive number $M$\ and for all $p\in S(A)$;

(2) $(\tint\limits_Gf(s)ds)a=\tint\limits_Gf(s)ads$\ for all $a$\ $\in A$;

(3) $\Phi (\tint\limits_Gf(s)ds)=\tint\limits_G\Phi \left( f(s)\right) ds$\
for any morphism of locally $C^{*}$-algebras $\Phi :A\rightarrow B$;

(4) $(\tint\limits_Gf(s)ds)^{*}=\tint\limits_Gf(s)^{*}ds$.
\end{quote}
\end{Lemma}

\begin{quote}
\smallskip
\end{quote}

\proof%
Let $p\in S(A)$. Then $\pi _p\circ f$ $\in C_c(G,A_p)$ and so there is a
unique element $\tint\limits_G(\pi _p\circ f)(s)ds$ in $A_p$ such that for
any non-degenerate $*$ -representation $(\varphi _p,H_{\varphi _p})$ of $A_p$
\[
\left\langle \varphi _p(\tint\limits_G(\pi _p\circ f)(s)ds)\xi ,\eta
\right\rangle =\tint\limits_G\left\langle \varphi _p((\pi _p\circ f)(s))\xi
,\eta \right\rangle ds 
\]
for all $\xi ,\eta $ in $H_{\varphi _p}$, see, for instance, Lemma 7 of [%
{\bf 11}].

To show that $(\tint\limits_G(\pi _p\circ f)(s)ds)_p$ is a coherent net in $%
A $, let $p,q\in S(A)$ with $p\geq q$. Then we have 
\begin{eqnarray*}
\pi _{pq}(\tint\limits_G(\pi _p\circ f)(s)ds) &=&\tint\limits_G\pi
_{pq}\left( (\pi _p\circ f)(s)\right) ds\; \\
&&\text{ussing Lemma 7 of [{\bf 11}]} \\
&=&\tint\limits_G(\pi _q\circ f)(s)ds.
\end{eqnarray*}
Therefore $(\tint\limits_G(\pi _p\circ f)(s)ds)_p\in A,$ and we define $%
\tint\limits_Gf(s)ds=(\tint\limits_G(\pi _p\circ f)(s)ds)_p$.

Suppose that there is another element $b$ in $A$ such that for any
non-degenerate $*$-representation $(\varphi ,H_{\varphi ^{}})$ of $A$ 
\[
\left\langle \varphi \left( b\right) \xi ,\eta \right\rangle
=\tint\limits_G\left\langle \varphi (f(s))\xi ,\eta \right\rangle ds 
\]
for all $\xi ,\eta $ in $H_{\varphi _{}}$.Then for any $p\in S(A)$ and for
any non-degenerate $*$-representation $(\varphi _p,H_{\varphi _p})$ of $A_p$ 
\[
\left\langle \varphi _p\left( \pi _p(b)\right) \xi ,\eta \right\rangle
=\tint\limits_G\left\langle \varphi _p((\pi _p\circ f)(s))\xi ,\eta
\right\rangle ds 
\]
for all $\xi ,\eta $ in $H_{\varphi _p}$. From these facts and Lemma 7 of [%
{\bf 11}], we conclude that 
\[
\pi _p(b)=\tint\limits_G(\pi _p\circ f)(s)ds 
\]
for all $p\in S(A)$. Therefore $b=\tint\limits_Gf(s)ds$ and the uniqueness
is proved.

Using Lemma 7 of [{\bf 11}] it is easy to check that $\tint\limits_Gf(s)ds$
satisfies the conditions $(1)-(4)$.%
\endproof%

\smallskip

Let $f,h$ in $C_c(G,A)$. It is easy to check that the map $(s,t)\rightarrow
f(t)\alpha _t\left( h(t^{-1}s)\right) $ from $G\times G$ to $A$ is an
element in $C_c(G\times G,A)$ and the relation 
\[
\left( f\times h\right) \left( s\right) =\tint\limits_Gf(t)\alpha _t\left(
h(t^{-1}s)\right) dt 
\]
defines an element in $C_c(G,A),$ called the convolution of $f$ and $h$.
Also it is not hard to check that $C_c(G,A)$ becomes a $*$-algebra with
convolution as product and involution defined by 
\[
f^{\sharp }(t)=\gamma (t)^{-1}\alpha _t\left( f(t^{-1})^{*}\right) 
\]
where $\gamma $ is the modular function on $G$.

For any $p\in S(A),$ define $N_p$ from $C_c(G,A)$ to $[0,\infty )$ by 
\[
N_p(f)=\tint\limits_Gp(f(s))ds. 
\]
Straightforward computations show that $N_p$, $p\in S(A)$ are
submultiplicative $*$-semi-norms on $C_c(G,A)$.

Let $L^1(G,A,\alpha )$ be the Hausdorff completion of $C_c(G,A)$ with
respect to the topology defined by the family of submultiplicative $*$%
-semi-norms $\{N_p\}_{p\in S(A)}$. Then by Theorem III 3.1 of [{\bf 7}] 
\[
L^1(G,A,\alpha )=\lim\limits_{\stackunder{p\in S(A)}{\leftarrow }}\left(
L^1(G,A,\alpha )\right) _p\; 
\]
where $\left( L^1(G,A,\alpha )\right) _p$ is the completion of the $*$
-algebra $C_c(G,A)/\ker \left( N_p\right) $ with respect to the norm $%
\left\| \cdot \right\| _p$ induced by $N_p$.

\smallskip

\begin{Lemma}
Let $(G,A,\alpha )$\ be a locally $C^{*}$-dynamical system such that $\alpha 
$\ is an inverse limit action. Then 
\[
\left( L^1(G,A,\alpha )\right) _p=L^1\left( G,A_p,\alpha ^{(p)}\right) 
\]
for all $p\in S(A)$, up to a topological algebraic $*$\ -isomorphism.
\end{Lemma}

\smallskip

\proof%
{\it \ }Let $p\in S(A)$ and $f$ in $C_c(G,A)$. Then 
\[
\left\| f+\ker (N_p)\right\| _p=\int\limits_Gp\left( f(s)\right)
ds=\int\limits_G\left\| \pi _p\left( f(s)\right) \right\| _pds=\left\| \pi
_p\circ f\right\| _1. 
\]
Therefore we can define a linear map $\psi _p$ from $C_c(G,A)/\ker \left(
N_p\right) $ to $C_c(G,A_p)$ by 
\[
\psi _p\left( f+\ker (N_p)\right) =\pi _p\circ f. 
\]

It is not hard to check that $\psi _p$ is a $*$ -morphism, and since $\psi
_p $ is an isometric $*$ -morphism from $C_c(G,A)/\ker (N_p)$ to $C_c(G,A_p)$%
, it can be uniquely extended to an isometric $*$ -morphism $\psi _p$ from $%
\left( L^1(G,A,\alpha )\right) _p$ to $L^1\left( G,A_p,\alpha ^{(p)}\right) $%
.

To show that $\psi _p$ is surjective, let $a\in A$ and $f\in C_c(G)$. Define 
$\widetilde{f}$ from $G$ to $A$ by $\widetilde{f}(s)=f(s)a$. Clearly $%
\widetilde{f}\in C_c(G,A)$ and 
\[
\psi _p\left( \widetilde{f}+\ker (N_p)\right) (s)=f(s)\pi _p(a) 
\]
for all $s$ in $G$. This implies that 
\[
A_p\otimes _{\text{alg }}C_c(G)\subseteq \psi _p\left( \left( L^1(G,A,\alpha
)\right) _p\right) \subseteq L^1(G,A_p,\alpha ^{(p)}) 
\]
whence, since $A_p\otimes _{\text{alg }}C_c(G)$ is dense in $%
L^1(G,A_p,\alpha ^{(p)})$ and since $\psi _p$ is an isometric $*$ -morphism,
we deduce that $\psi _p$ is surjective and the proposition is proved.%
\endproof%

\smallskip

\begin{Corollary}
Let $(G,A,\alpha )$\ be a locally $C^{*}$-dynamical system such that $\alpha 
$\ is an inverse limit action. Then 
\[
L^1(G,A,\alpha )=\lim\limits_{\stackrel{\leftarrow }{p\in S(A)}}L^1\left(
G,A_p,\alpha ^{(p)}\right) 
\]
up to an algebraic and topological $*$-isomorphism.
\end{Corollary}

\smallskip

\begin{Remark}
{\em If }$\{e_i\}_{i\in I}${\em \ is an approximate unit for }$A${\em \ and }%
$\{f_j\}_{j\in J}${\em \ is an approximate unit for }$L^1(G)${\em , then }$\{%
\widetilde{f}_{(i,j)}\}_{(i,j)\in I\times J}${\em , where }$\widetilde{f}%
_{(i,j)}(s)=f_j(s)e_i,${\em \ }$s\in G${\em , is an approximate unit for }$%
L^1(G,A,\alpha ),${\em \ see Lemma XIV.1.2 of [{\bf 7}]. Then by Definition
5.1 of [{\bf 1}], }${\em we}$ {\em can construct the enveloping algebra of} $%
L^1(G,A,\alpha ).$
\end{Remark}

\smallskip

\begin{Definition}
{\em A covariant representation of }$(G,A,\alpha )${\em \ is a triple }$%
(\varphi ,u,H)${\em , where }$(\varphi ,H)${\em \ is a }$*${\em %
-representation of }$A${\em \ and }$\left( u,H\right) ${\em \ is a unitary
representation of }$G${\em \ such that } 
\[
\varphi (\alpha _t(a))=u_t\varphi (a)u_t^{*} 
\]
{\em for all }$t\in G${\em \ and for all }$a\in A${\em .}

{\em We say that the covariant representation }$(\varphi ,u,H)${\em \ of }$%
(G,A,\alpha )${\em \ is non-degenerate if the }$*${\em -representation }$%
(\varphi ,H)${\em \ of }$A${\em \ is non-degenerate.}
\end{Definition}

\begin{Remark}
$\smallskip (1).{\em \ }${\em If }$(\varphi ,u,H)${\em \ is a covariant
representation of }$(G,A,\alpha )${\em \ such that }$\left\| \varphi
(a)\right\| \leq p(a)$ {\em \ for all }$a\in A${\em , then there is a unique
covariant representation }$(\varphi _p,u,H)${\em \ of the }$C^{*}${\em %
-dynamical system }$\left( G,A_p,\alpha ^{(p)}\right) ${\em \ such that }$%
\varphi _p\circ \pi _p=\varphi ${\em .}

$(2).$ {\em If }$(\varphi _p,u,H)${\em \ is a covariant representation of
the }$C^{*}${\em -dynamical system }$\left( G,A_p,\alpha ^{(p)}\right) ${\em %
, then }$(\varphi _p\circ \pi _p,u,H)${\em \ is a covariant representation
of the locally }$C^{*}${\em -dynamical system }$\left( G,A,\alpha \right) $%
{\em .}
\end{Remark}

\smallskip

If $R(G,A,\alpha )$ denotes the non-degenerate covariant representations of $%
(G,A,\alpha ),$ then it is easy to check that 
\[
R(G,A,\alpha )=\bigcup\limits_{p\in S(A)}R_p(G,A,\alpha ) 
\]
where $R_p(G,A,\alpha )=\left\{ (\varphi ,u,H)\in R(G,A,\alpha );\left\|
\varphi (a)\right\| \leq p(a)\text{ for all }a\in A\right\} $. Also it is
easy to check that the map $\varphi _p\mapsto \varphi _p\circ \pi _p$ from $%
R(G,A_p,\alpha ^{(p)})$ to $R_p(G,A,\alpha )$ is bijective.

\begin{Proposition}
Let $(G,A,\alpha )$\ be a locally $C^{*}$-dynamical system such that $\alpha 
$\ is an inverse limit action. Then there is a bijection between the
covariant non-degenerate representations of $(G,A,\alpha )$\ and the
non-degenerate $*$-representations of $L^1(G,A,\alpha ).$
\end{Proposition}

\smallskip

\proof%
{\it \ }Let $(\varphi ,u,H)\in R(G,A,\alpha )$. Then, there is $p\in S(A)$
and $(\varphi _p,u,H)\in R(G,A_p,\alpha ^{(p)})$ such that $\varphi =\varphi
_p\circ \pi _p$. Since $(\varphi _p,u,H)\in R(G,A_p,\alpha ^{(p)})$ there is
a unique non-degenerate $*$-representation $(\varphi _p\times u,H)$ of $%
L^1(G,A_p,\alpha ^{(p)})$ such that 
\[
(\varphi _p\times u)\left( f\right) =\int\limits_G\varphi _p\left(
f(t)\right) u_tdt\; 
\]
for all $f\in L^1\left( G,A_p,\alpha ^{(p)}\right) $, see, for instance,
Proposition 7.6.4 of [{\bf 8}].

Let $\varphi \times u=(\varphi _p\times u)\circ $ $\widetilde{\pi }_p$,
where $\widetilde{\pi }_p$ is the canonical map from $L^1(G,A,\alpha )$ to $%
L^1(G,A_p,\alpha ^{(p)}),$ $\widetilde{\pi }_p(f)=\pi _p\circ f$ for all $f$
in $L^1(G,A,\alpha )$. Then, clearly $(\varphi \times u,H)$ is a
non-degenerate $*$-representation of $L^1(G,A,\alpha )$ and moreover, 
\[
(\varphi \times u)\left( f\right) =(\varphi _p\times u)(\pi _p\circ
f)=\int\limits_G\varphi _p\left( (\pi _p\circ f)(t)\right)
u_tdt=\int\limits_G\varphi (f(t))u_tdt 
\]
for all $f\in L^1(G,A,\alpha )$. Thus we have obtained a map $(\varphi
,u,H)\rightarrow (\varphi \times u,H)$ from $R(G,A,\alpha )$ to $%
R(L^1(G,A,\alpha ))$. To show that this map is bijective, let $(\Phi ,H)$ be
a non-degenerate $*$ -representation of $L^1(G,A,\alpha )$. Then there is $%
p\in S(A)$ and a non-degenerate $*$-representation $(\Phi _p,H)$ of $%
L^1(G,A_p,\alpha ^{(p)})$ such that $\Phi =\Phi _p\circ \pi _p$. By
Proposition 7.6.4 of [{\bf 8}] there is a unique non-degenerate covariant
representation $(\varphi _p,u,H)$ of $(G,A_p,\alpha ^{(p)})$ such that $%
(\phi _p,H)=(\varphi _p\times u,H)$. Therefore there is a non-degenerate
covariant representation $(\varphi ,u,H)\ $of $(G,A,\alpha )$, where $%
\varphi =\varphi _p\circ \pi _p$, such that $(\Phi ,H)=(\varphi \times u,H)$.

To show that $(\varphi ,u,H)$ is unique, let $(\psi ,v,K)$ be another
non-degenerate covariant representation of $\left( G,\alpha ,A\right) $ such
that $(\psi \times v,K)=(\Phi ,H)$. Then there is $q\in S(A)$ with $q\geq p$
such that $(\psi ,v,K)\in R_q\left( G,A,\alpha \right) $ and $\left( \Phi
,K\right) \in R_q\left( L^1(G,A,\alpha )\right) $. Therefore $\Phi =\Phi
_q\circ \widetilde{\pi }_q$ with $(\Phi _q,H)\in R\left( L^1(G,A_q,\alpha
^{(q)})\right) $ and $\psi =\psi _q\circ \pi _q$ with $(\psi _q,v,K)\in
R\left( G,A_q,\alpha ^{(q)}\right) $ and moreover, $(\Phi _q,H)=(\psi
_q\times v,K)$.

On the other hand, $(\varphi _p\circ \pi _{pq},u,H)\in $ $R\left(
G,A_q,\alpha ^{(q)}\right) $ and 
\begin{eqnarray*}
((\varphi _p\circ \pi _{pq})\times u)(f) &=&\int\limits_G\left( \varphi
_p\circ \pi _{pq}\right) (f(t))u_tdt=\int\limits_G\varphi _p\left( 
\widetilde{\pi }_{pq}(f)(t)\right) u_tdt \\
&=&\phi _p\left( \widetilde{\pi }_{pq}(f)\right) =\left( \phi _p\circ 
\widetilde{\pi }_{pq}\right) (f)=\phi _q(f)
\end{eqnarray*}
for all $f\in L^1(G,A_q,\alpha ^{(q)})$. From these facts and Proposition
7.6.4 of [{\bf 8}], we conclude that the covariant representations $\left(
\psi _q,v,K\right) $ and $\left( \varphi _p\circ \pi _{pq},u,H\right) $ of $%
\left( G,A_q,\alpha ^{(q)}\right) $ coincide, and so the covariant
representations $\left( \psi ,v,K\right) $ and $\left( \varphi ,u,H\right) $
of $\left( G,A,\alpha \right) $ coincide.%
\endproof%

\smallskip

\begin{Definition}
{\em Let }$\left( G,A,\alpha \right) ${\em \ be a locally }$C^{*}${\em %
-dynamical system such that }$\alpha ${\em \ is an inverse limit action. The
crossed product of }$A${\em \ by the action }$\alpha ,$ {\em denoted by }$%
A\times _{\alpha _{}}G,${\em \ is the enveloping algebra of the complete
locally }$m${\em \ -convex }$*${\em \ -algebra }$L^1(G,A,\alpha ).$
\end{Definition}

\smallskip

\begin{Remark}
{\em By Corollary 3.9 and Corollary 5.3 of [{\bf 2}], }$A\times _{\alpha
_{}}G${\em \ is a locally }$C^{*}${\em -algebra and } 
\[
A\times _{\alpha _{}}G=\lim\limits_{\stackrel{\longleftarrow }{p\in S(A)}%
}A_p\times _{\alpha ^{(p)}}G 
\]
{\em up to an isomorphism of locally }$C^{*}${\em -algebras.}
\end{Remark}

{\em \smallskip }

\begin{Proposition}
Let $\left( G,A,\alpha \right) $\ be a locally $C^{*}$-dynamical system such
that $\alpha $\ is an inverse limit action. Then there is a bijection
between non-degenerate covariant representations of $\left( G,A,\alpha
\right) $\ and the non-degenerate representations of $A\times _{\alpha _{}}G$%
.
\end{Proposition}

\smallskip

\proof%
{\it \ }Since $A\times _{\alpha _{}}G$ is the enveloping locally $C^{*}$%
-algebra of the complete locally $m$ -convex $*$ -algebra $L^1(G,A,\alpha ),$
there is a bijection between the non-degenerate representations of $A\times
_{\alpha _{}}G$ and the non-degenerate representations of $L^1(G,A,\alpha ),$
[{\bf 2, }pp. 37]. From this fact and Proposition 3.13 we conclude that
there is a bijection between the non-degenerate representations of $A\times
_{\alpha _{}}G$ and the non-degenerate covariant representations of $\left(
G,A,\alpha \right) .$%
\endproof%

\smallskip

For each $p\in S(A)$, we denote by $\left( \varphi _{p,u},H_{p,u}\right) $
the universal representation of $A_p$ and by $\left( \varphi
_p,H_{p,u}\right) $ the representation of $A$ associated with $\left(
\varphi _{p,u},H_{p,u}\right) $ (that is, $\varphi _p=\varphi _{p,u}\circ
\pi _p$).

\begin{Lemma}
Let $\left( G,A,\alpha \right) $\ be a locally $C^{*}$-dynamical system such
that $\alpha $\ is an inverse limit action. Then $\left( \widetilde{\varphi
_p},\lambda ,L^2\left( G,H_{p,u}\right) \right) $, where 
\[
\widetilde{\varphi _p}\left( a\right) \left( \xi \right) \left( t\right)
=\varphi _p\left( \alpha _{t^{-1}}\left( a\right) \right) \left( \xi \left(
t\right) \right) 
\]
and 
\[
\lambda _s\left( \xi \right) \left( t\right) =\xi \left( s^{-1}t\right) 
\]
for all $a$\ in $A$, $\xi $\ in $L^2\left( G,H_{p,u}\right) $\ and $s,t$\ in 
$G,$\ is a non-degenerate covariant representation of $\left( G,A,\alpha
\right) .$
\end{Lemma}

\smallskip

\proof%
{\it \ }It is a simple verification.%
\endproof%

\smallskip

Let $p\in S(A)$. The map $r_p:L^1\left( G,A,\alpha \right) \rightarrow
[0,\infty )$ defined by 
\[
r_p(f)=\left\| \left( \widetilde{\varphi _p}\times \lambda \right) \left(
f\right) \right\| 
\]
is a $C^{*}$-semi-norm on $L^1\left( G,A,\alpha \right) $ with the property
that $r_p(f)\leq N_p(f)$ for all $f$ in $L^1\left( G,A,\alpha \right) $.

Let $I=\tbigcap\limits_{p\in S(A)}\ker \left( r_p\right) $. Clearly $I$ is a
closed two-sided ideal of $L^1\left( G,A,\alpha \right) $ and $L^1\left(
G,A,\alpha \right) /I$ is a pre-locally $C^{*}$-algebra with respect to the
topology determined by the family of $C^{*}$-semi-norms $\{\widehat{r}%
_p\}_{p\in S(A)},$ $\widehat{r}_p(f+I)=\inf \{r_p(f+h);h\in I\}$.

\smallskip

\begin{Definition}
{\em The reduced crossed product of }$A${\em \ by the action }$\alpha ,${\em %
\ denoted by }$A\times _{\alpha ,r}G${\em , is the Hausdorff completion of }$%
\left( L^1\left( G,A,\alpha \right) ,\{r_p\}_{p\in S(A)}\right) ${\em \ (
that is, }$A\times _{\alpha ,r}G${\em \ is the completion of the pre-locally 
}$C^{*}${\em -algebra }$\left( L^1\left( G,A,\alpha \right) /I,\{\widehat{r}%
_p\}_{p\in S(A)}\right) ${\em ).}
\end{Definition}

\smallskip

\begin{Lemma}
{\it \ }Let $\left( G,A,\alpha \right) $\ be a locally $C^{*}$-dynamical
system such that $\alpha $\ is an inverse limit action. Then 
\[
\left( A\times _{\alpha ,r}G\right) _p=A_p\times _{\alpha ^{(p)},r}G 
\]
for all $p\in S(A)$, up to an isomorphism of $C^{*}$-algebras.
\end{Lemma}

\smallskip

\proof%
Let $p\in S(A)$. If $f\in L^1\left( G,A,\alpha \right) $, then we have 
\begin{eqnarray*}
\left\| \left( f+I\right) +\ker (\widehat{r}_p)\right\| _{\widehat{r}_p} &=&%
\widehat{r}_p(f+I)=\inf \{\left\| \left( \widetilde{\varphi _p}\times
\lambda \right) \left( f+h\right) \right\| ;h\in I\} \\
&=&\inf \{\left\| \left( \widetilde{\varphi _p}\times \lambda \right) \left(
f\right) \right\| ;h\in I\}=r_p(f)=\left\| f+\ker (r_p)\right\| _{r_p}.
\end{eqnarray*}
From this relation, we conclude that $\left( A\times _{\alpha ,r}G\right) _p$
is isomorphic to the completion of $L^1\left( G,A,\alpha \right) /\ker (r_p)$
with respect to the $C^{*}$-norm induced by $r_p$.

On the other hand, $A_p\times _{\alpha ^{(p)},r}G$ is the completion of $%
L^1\left( G,A_p,\alpha ^{(p)}\right) /I_p,$ where $I_p=\{f\in L^1\left(
G,A_p,\alpha ^{(p)}\right) /\left( \widetilde{\varphi _{p,u}}\times \lambda
\right) \left( f\right) =0\},$ with respect to the norm $\left\| \cdot
\right\| $ $^{^{\prime }}$ given by $\left\| f+I_p\right\| ^{\prime
}=\left\| \left( \widetilde{\varphi _{p,u}}\times \lambda \right) \left(
f\right) \right\| \leq \left\| f\right\| _1$. But the completion of $%
L^1\left( G,A,\alpha \right) $\ $/\ker (r_p)$ with respect to the norm $%
\left\| \cdot \right\| _{r_p}$ is isomorphic to the completion of $L^1\left(
G,A,\alpha ^{(p)}\right) /I_p$ with respect to the norm $\left\| \cdot
\right\| ^{\prime },$ since 
\begin{eqnarray*}
\left\| f+\ker (r_p)\right\| _{r_p} &=&r_p(f)=\left\| \left( \widetilde{%
\varphi _p}\times \lambda \right) \left( f\right) \right\| \\
&=&\left\| \left( \widetilde{\varphi _{p,u}}\times \lambda \right) \left(
\pi _p\circ f\right) \right\| =\left\| \widetilde{\pi }_p(f)+I_p\right\|
^{\prime }
\end{eqnarray*}
for all $f\in L^1(G,A,\alpha )$. Therefore the $C^{*}$ -algebras $\left(
A\times _{\alpha ,r}G\right) _p$ and $A_p\times _{\alpha ^{(p)},r}G$ are
isomorphic.%
\endproof%

\smallskip

\begin{Corollary}
If $\left( G,A,\alpha \right) $\ is a locally $C^{*}$-dynamical system such
that $\alpha $\ is an inverse limit action then 
\[
A\times _{\alpha ,r}G=\lim\limits_{\stackunder{p\in S(A)}{\leftarrow }%
}A_p\times _{\alpha ^{(p)},r}G 
\]
up to an isomorphism of locally $C^{*}$-algebras.
\end{Corollary}

\smallskip

\section{ Actions of a Kac $C^{*}$-algebra on a locally $C^{*}$-algebra}

Let $C^{*}{\bf K}$ $=(B,d,j,\varphi )$ be a Kac $C^{*}$-algebra and let $A$
be a locally $C^{*}$-algebra.

\begin{Definition}
{\em An action of }$C^{*}{\bf K}${\em \ on }$A${\em \ is an injective }$S$%
{\em \ -morphism }$\alpha $ {\em \ from }$A${\em \ to }$M(A,B)${\em \ such
that } 
\[
\left( \overline{\alpha \otimes \text{id}_B}\right) \circ \alpha =\left( 
\overline{\text{id}_A\otimes (\sigma _B\circ d)}\right) \circ \alpha . 
\]
{\em An action }$\alpha ${\em \ of }$C^{*}{\bf K}${\em \ on }$A${\em \ is an
inverse limit action if we can write }$A${\em \ as an inverse limit }$%
\lim\limits_{\stackunder{\delta \in \Delta }{\leftarrow }}A_{\delta ^{}}$%
{\em \ of }$C^{*}${\em -algebras such a way that there are actions }$\alpha
^{\left( \delta \right) }${\em \ of }$C^{*}{\bf K}${\em \ on }$A_{\delta
^{}},${\em \ }$\delta \in \Delta $ {\em \ such that }$\alpha =${\em \ }$%
\lim\limits_{\stackunder{\delta \in \Delta }{\leftarrow }}\alpha ^{\left(
\delta \right) }.$

{\em Two actions }$\alpha _1${\em \ and }$\alpha _2${\em \ of }$C^{*}{\bf K}$%
{\em \ on the locally }$C^{*}${\em -algebras }$A_1${\em \ respectively }$A_2$%
{\em \ are said to be equivalent if there is an isomorphism of locally }$%
C^{*}${\em -algebras }$\Phi :A_1\rightarrow A_2${\em \ such that }$\alpha
_2\circ \Phi =\left( \overline{\Phi \otimes \text{id}_B}\right) \circ \alpha
_1.$
\end{Definition}

\smallskip

\begin{Proposition}
Let $G$\ be a locally compact group. If $\alpha $\ is an action of $C^{*}%
{\bf K}_G^a$\ on $A$,\ then the map $\Sigma \left( \alpha \right) $\ that
applies $t\in G$ to a map $\Sigma \left( \alpha \right) _t$\ from $A$\ to $A$%
\ defined by $\Sigma \left( \alpha \right) _t\left( a\right) =\alpha
(a)\left( t^{-1}\right) $,\ is a continuous action of $G$\ on $A.$
\end{Proposition}

\smallskip

\proof%
Since $\alpha $ is a continuous $*$-morphism from $A$ to $C_b(G,A)$, $\Sigma
\left( \alpha \right) _t$ is a continuous $*$-morphism from $A$ to $A$ for
each $t$ in $G$. Using the same arguments as in the proof of Proposition
5.1.5 of [{\bf 13}], it is not difficult to see that $\Sigma \left( \alpha
\right) _t$ is invertible and moreover, $\left( \Sigma \left( \alpha \right)
_t\right) ^{-1}=\Sigma \left( \alpha \right) _{t^{-1}}$ for all $t$ in $G$.
Therefore $\Sigma \left( \alpha \right) _t\in $Aut$(A)$ for each $t$ in $G$.

To show that the map $\left( t,a\right) \rightarrow $ $\Sigma \left( \alpha
\right) _t\left( a\right) $ from $G\times A$ to $A$ is continuous, let $%
\left( t_0,a_0\right) \in G\times A$ and let $W_{p,\varepsilon }=\{a\in
A;p(a-$ $\Sigma \left( \alpha \right) _{t_0}\left( a_0\right) )<\varepsilon
\}$ be a neighborhood of $\Sigma \left( \alpha \right) _{t_0}\left(
a_0\right) $ . Since $\alpha \left( a_0\right) \in C_b\left( G,A\right) $,
there is a neighborhood $U_0$ of $t_0$ such that 
\[
p\left( \alpha \left( a_0\right) (t^{-1})-\alpha \left( a_0\right)
(t_0^{-1})\right) <\frac \varepsilon 2 
\]
for all $t$ in $U_0$, and since $\alpha $ is a continuous $*$-morphism,
there is a neighborhood $V_0$ of $a_0$ such that 
\[
\left\| \alpha \left( a\right) -\alpha \left( a_0\right) \right\| _p=\sup
\{p(\alpha \left( a\right) (t)-\alpha \left( a_0\right) (t));t\in G\}<\frac
\varepsilon 2 
\]
for all $a$ in $V_0$. Then 
\begin{eqnarray*}
p\left( \Sigma \left( \alpha \right) _t\left( a\right) -\Sigma \left( \alpha
\right) _{t_0}\left( a_0\right) \right) &\leq &p\left( \alpha \left(
a\right) \left( t^{-1}\right) -\alpha \left( a_0\right) \left( t^{-1}\right)
\right) \\
&&+p\left( \alpha \left( a_0\right) \left( t^{-1}\right) -\alpha \left(
a_0\right) (t_0^{-1})\right) \\
&\leq &\left\| \alpha \left( a\right) -\alpha \left( a_0\right) \right\|
_p+\frac \varepsilon 2<\varepsilon
\end{eqnarray*}
for all $(t,a)\in U_0\times V_0$ and the proposition is proved. 
\endproof%

\smallskip

\begin{Remark}
{\em According to Proposition 4.2, we can define a map }$\Sigma ${\em \ from
the set of all actions of }$C^{*}{\bf K}_G^a${\em \ on }$A${\em \ to the set
of all continuous actions of }$G${\em \ on }$A${\em \ by }$\alpha
\rightarrow \Sigma \left( \alpha \right) ${\em . Moreover, }$\Sigma ${\em \
is injective.}
\end{Remark}

{\em \smallskip }

The following proposition is a generalization of Proposition 5.1.5 of [{\bf %
13}] for inverse limit actions of locally compact groups on locally $C^{*}$%
-algebras.

\begin{Proposition}
Let $G$\ be a locally compact group. Then the map $\Sigma $\ defined in
Proposition 4.2 is a bijective correspondence between the set of all inverse
limit actions of $C^{*}{\bf K}_G^a$\ on $A$\ and the set of all continuous
inverse limit actions of $G$\ on $A.$
\end{Proposition}

\smallskip

\proof%
Let $\alpha $ be an inverse limit action of $C^{*}{\bf K}_G^a$ on $A$. Then $%
A$ may be written as an inverse limit $\lim\limits_{\stackunder{\delta \in
\Delta }{\leftarrow }}A_{\delta ^{}}$ of $C^{*}$-algebras and there are
actions $\alpha ^{\left( \delta \right) }$ of $C^{*}{\bf K}_G^a$ on $%
A_{\delta _{}},$ $\delta \in \Delta $ such that $\alpha =\lim\limits_{%
\stackunder{\delta \in \Delta }{\leftarrow }}\alpha ^{\left( \delta \right)
} $.

According to Proposition 5.1.5 of [{\bf 13}], for each $\delta \in \Delta $
there is a continuous action $\Sigma (\alpha ^{\left( \delta \right) })$ of $%
G$ on $A_{\delta ^{}}$ such that $\Sigma (\alpha ^{\left( \delta \right)
})_t(a_{\delta ^{}})=\alpha ^{\left( \delta \right) }(a_{\delta
_{}})(t^{-1}) $ for all $a_{\delta ^{}}$ in $A_{\delta _{}}$ and for all $t$
in $G$. Since $\{\alpha ^{\left( \delta \right) }\}_{\delta \in \Delta }$ is
an inverse system of morphisms of $C^{*}$-algebras, it is not difficult to
check that $\{\Sigma (\alpha ^{\left( \delta \right) })_t\}_{\delta \in
\Delta }$ is an inverse system of $C^{*}$-isomorphisms for each $t$ in $G$.
Also it is easy to check that $\Sigma \left( \alpha \right) _t=\lim\limits_{%
\stackunder{\delta \in \Delta }{\leftarrow }}$ $\Sigma (\alpha ^{\left(
\delta \right) })_t$ for each $t$ in $G$.

To show that $\Sigma $ is surjective, let $\beta $ be a continuous inverse
limit action of $G$ on $A$. Then $A$ may be written as an inverse limit $%
A=\lim\limits_{\stackunder{\delta \in \Delta }{\leftarrow }}A_{\delta ^{}}$
of $C^{*}$-algebras and there are continuous actions $\beta ^{\left( \delta
\right) }$ of $G$ on $A_{\delta ^{}}$ such that $\beta _t=\lim\limits_{%
\stackunder{\delta \in \Delta }{\leftarrow }}$ $\beta _t^{\left( \delta
\right) }$ for each $t$ in $G$. By Proposition 5.1.5 of [{\bf 13}], for each 
$\delta \in \Delta $ there is an action $\alpha ^{\left( \delta \right) }$
of $C^{*}{\bf K}_G^a$ on $A_{\delta ^{}}$ such that $\Sigma (\alpha ^{\left(
\delta \right) })=\beta ^{\left( \delta \right) }$. It is not difficult to
verify that $\{\alpha ^{\left( \delta \right) }\}_{\delta \in \Delta }$ is
an inverse system of injective $S$-morphisms of $C^{*}$-algebras. Let $%
\alpha =\lim\limits_{\stackunder{\delta \in \Delta }{\leftarrow }}$ $\alpha
^{\left( \delta \right) }$. Then $\alpha $ is an injective $S$-morphism of
locally $C^{*}$-algebras and 
\begin{eqnarray*}
\left( \overline{\alpha \otimes \text{id}_{C_0(G)}}\right) \circ \alpha
&=&\lim\limits_{\stackunder{\delta \in \Delta }{\leftarrow }}\left( 
\overline{\alpha ^{\left( \delta \right) }\otimes \text{id}_{C_0(G)}}\right)
\circ \alpha ^{\left( \delta \right) } \\
&=&\lim\limits_{\stackunder{\delta \in \Delta }{\leftarrow }}\left( 
\overline{\text{id}_{A_p}\otimes \sigma _{C_0(G)}\circ d_G^a}\right) \circ
\alpha ^{\left( \delta \right) } \\
&=&\left( \overline{\text{id}_A\otimes \sigma _{C_0(G)}\circ d_G^a}\right)
\circ \alpha .
\end{eqnarray*}
Therefore $\alpha $ is an inverse limit action of $C^{*}{\bf K}_G^a$ on $A$
and $\Sigma \left( \alpha \right) =\beta $. Thus we showed that $\Sigma $ is
bijective.%
\endproof%

\smallskip

\begin{Corollary}
{\it \ }If $G$\ is compact, then any action of $C^{*}{\bf K}_G^a$\ on $A$\
is an inverse limit action.
\end{Corollary}

\smallskip

\proof%
Let $\alpha $ be an action of $C^{*}{\bf K}_G^a$ on $A$. By Proposition 4.2, 
$\Sigma \left( \alpha \right) $ is a continuous action of $G$ on $A$ which
is an limit inverse action, since the group $G$ is compact, Lemma 3.6. From
this fact and Proposition 4.4 we conclude that $\alpha $ is an inverse limit
action.%
\endproof%

\smallskip

\section{The Takai duality theorem}

Let $G$ be a locally compact group and let $A$ be a locally $C^{*}$-algebra.

\begin{Lemma}
Let $\alpha $\ be an inverse limit action of $G$\ on $A$. Then the reduced
crossed product of $A$\ by the action $\alpha $\ is isomorphic to the
locally $C^{*}$-subalgebra of $M\left( A\otimes {\cal L}\left( L^2\left(
G\right) \right) \right) $\ generated by $\{\alpha (a)\left( 1_{M(A)}\otimes
\lambda (f)\right) ;a\in A,f\in C_c\left( G\right) \}$, where $\lambda $\ is
the left regular representation of $L^1(G)$.
\end{Lemma}

\smallskip

\proof%
Let $p\in S(A)$. By Remark 5.2.1.1 of [{\bf 13}], the map $\Phi _p$ from the 
$C^{*}$-subalgebra of $M\left( A_p\otimes {\cal L}\left( L^2\left( G\right)
\right) \right) $ generated by $\{\alpha ^{(p)}(a_p)\left( 1_{M(A_p)}\otimes
\lambda (f)\right) ;a_p\in A_p,f\in C_c\left( G\right) \}$ to $A_p\times
_{\alpha ^{(p)},r}G$, that applies $\alpha ^{(p)}(a_p)(1_{M(A_p)}\otimes
\lambda (f))$ to $\widetilde{f}+I_p,$ where $\widetilde{f}(t)=f(t)a_p,$ $%
t\in G,$ see the proof of Lemma 3.19, is an isomorphism of $C^{*}$-algebras.

If $\pi _{pq}^{\prime },$ $p,$ $q\in S(A),$ $p\geq q$ are the connecting
maps\ of\ the\ inverse system\ \ \ $\{M(A_p\otimes {\cal L}\left( L^2\left(
G\right) \right) )\}_{p\in S(A)}$ and $\widehat{\pi }_{pq},$ $p,q\in S(A),$ $%
p\geq q$ are the connecting maps of the inverse system $\{A_p\times _{\alpha
^{(p)},r}G\}_{p\in S(A)}$, then we have 
\begin{eqnarray*}
(\Phi _q\circ \pi _{pq}^{\prime })(\alpha ^{(p)}(a_p)(1_{M(A_p)}\otimes
\lambda (f))) &=&\Phi _q(\alpha ^{(q)}(\pi _{pq}(a_p))(1_{M(A_q)}\otimes
\lambda (f))) \\
&=&\pi _{pq}(a_p)\otimes f+I_p=\widetilde{\pi }_{pq}(a_p\otimes f)+I_p \\
&=&\widehat{\pi }_{pq}(a_p\otimes f+I_p) \\
&=&(\widehat{\pi }_{pq}\circ \Phi _p)(\alpha ^{(p)}(a_p)(1_{M(A_p)}\otimes
\lambda (f)))
\end{eqnarray*}
for all $a_p$ in $A_p$, for all $f$ in $C_c(G)$ and for all $p,q\in S(A)$
with $p\geq q$. Therefore $\{\Phi _p\}_{p\in S(A)}$ is an inverse system of
isomorphisms of $C^{*}$-algebras and the lemma is proved.%
\endproof%

\smallskip

\begin{Definition}
{\em A coaction of }$G${\em \ on }$A${\em \ is an action }$\beta ${\em \ of }%
$C^{*}{\bf K}_G^s${\em \ on }$A.$ {\em We say that a coaction }$\beta ${\em %
\ of }$G${\em \ on }$A${\em \ is an inverse limit coaction if it\ is an
inverse limit action of }$C^{*}{\bf K}_G^s${\em \ on }$A.$

{\em The reduced crossed product of }$A${\em \ by the coaction }$\beta ,$%
{\em \ denoted by }$A\times _{\beta _{}}G,${\em \ is the locally }$C^{*}$%
{\em -subalgebra of }$M(A\otimes ${\em \ }${\cal L}\left( L^2\left( G\right)
\right) )${\em \ generated by }$\{\beta (a)\left( 1_{M(A)}\otimes f\right)
;a\in A,f\in C_c\left( G\right) \}.$
\end{Definition}

\smallskip

\begin{Remark}
{\em Let }$\beta =\lim\limits_{\stackunder{\delta \in \Delta }{\leftarrow }%
}\beta ^{(\delta )}${\em \ be an inverse limit coaction of }$G${\em \ on }$A$%
{\em \ such that the connecting maps of the inverse system }$\{A_{\delta
^{}}\}_{\delta \in \Delta }${\em \ are all surjective. Then, by Theorem 3.14
of [{\bf 10} ] }

\[
M(A\otimes {\em \ }{\cal L}\left( L^2\left( G\right) \right) )=\lim\limits_{%
\stackunder{\delta \in \Delta }{\leftarrow }}M(A_{\delta _{}}\otimes {\em \ }%
{\cal L}\left( L^2\left( G\right) \right) ) 
\]
{\em up to an isomorphism of locally }$C^{*}${\em -algebras, and} {\em by
Lemma III 3.2 of [{\bf 7}], } 
\[
A\times _{\beta _{}}G=\lim\limits_{\stackunder{\delta \in \Delta }{%
\leftarrow }}A_{\delta _{}}\times _{\beta ^{(\delta )}}G 
\]
{\em up to an isomorphism of locally }$C^{*}${\em -algebras.}
\end{Remark}

\smallskip

\begin{Remark}
{\em Let }$G${\em \ be a commutative locally compact group. Exactly as in
the proof of Proposition 5.1.6 of [{\bf 13}] we show that if }$\beta ${\em \
is an inverse limit coaction of }$G${\em \ on }$A${\em , then }$\beta
^{\prime }=(${\em id}$_A\otimes ${\em ad}${\cal F})\circ \beta ${\em , where 
}${\cal F}${\em \ is the Fourier-Plancherel isomorphism from }$L^2(G)${\em \
onto }$L^2(\widehat{G})${\em , is an inverse limit action of }$\widehat{G}$%
{\em \ on }$A${\em \ and conversely, if }$\alpha ${\em \ is an inverse limit
action of }$\widehat{G}${\em \ on }$A${\em \ then }$\alpha ^{\prime }=(${\em %
id}$_A\otimes ${\em ad}${\cal F}^{*})\circ \alpha $ {\em \ is an inverse
limit coaction of }$G${\em \ on }$A${\em . Therefore an inverse limit
coaction of }$G$ {\em can be identified with an inverse limit action of }$%
\widehat{G}${\em \ and id}$_A\otimes ${\em ad}${\cal F}${\em \ is an
isomorphism between }$A\times _{\beta _{}}G${\em \ and }$A\times _{\beta
^{\prime },r}\widehat{G}${\em .}
\end{Remark}

\smallskip

The following proposition is a generalization of Theorem 5.2.6 of [{\bf 13}]
for inverse limit actions of a locally compact group on a locally $C^{*}$%
-algebra.

\begin{Proposition}
Let $A$ be a locally $C^{*}$-algebra and let $G$ be a locally compact group.

$(1).$ If $\alpha $\ is an inverse limit action of $G$\ on $A$,\ then there
is an inverse limit coaction $\widehat{\alpha }$\ of $G$\ on $A\times
_{\alpha ,r}G$, called the dual coaction associated to $\alpha $, such that 
\begin{equation}
\widehat{\alpha }(\alpha (a)(1_{M(A)}\otimes \lambda (f)))=(\alpha
(a)\otimes 1_G)(1_{M(A)}\otimes d_G^s(\lambda (f)))  \tag{*}
\end{equation}
for all $a$\ in $A$\ and for all $f$\ in $C_c(G)$.

$(2).$ If $\beta =\lim\limits_{\stackunder{\delta \in \Delta }{\leftarrow }%
}\beta ^{(\delta )}$\ is an inverse limit coaction of $G$\ on $A$ such that
the connecting maps of the inverse system $\{A_{\delta ^{}}\}_{\delta \in
\Delta }$\ are all surjective, then there is an inverse limit action $%
\widehat{\beta }$\ of $G$\ on $A\times _{\beta _{}}G$, called the dual
action associated to $\beta $, such that 
\begin{equation}
\widehat{\beta }(\beta (a)(1_{M(A)}\otimes f))=(\beta (a)\otimes
1_G)(1_{M(A)}\otimes (\overline{\text{id}_{C_0(G)}\otimes j_G^a})d_G^a(f)) 
\tag{**}
\end{equation}
for all $a$\ in $A$\ and for all $f$\ in $C_c(G)$.
\end{Proposition}

\smallskip

\proof%
$(1).$ Since $\alpha $\ is an inverse limit action, $\alpha =$ $\lim\limits_{%
\stackunder{p\in S(A)}{\leftarrow }}\alpha ^{(p)}$, where $\alpha ^{(p)}$ is
a continuous action of $G$ on $A_p$. By Theorem 5.2.6 (i) of [{\bf 13}], for
each $p\in S(A)$ there is a dual coaction $\widehat{\alpha }^{(p)}$ of $G$
on $A_p\times _{\alpha ^{(p)},r}G$ such that 
\[
\widehat{\alpha }^{(p)}(\alpha ^{(p)}(a_p)(1_{M(A_p)}\otimes \lambda
(f)))=(\alpha ^{(p)}(a_p)\otimes 1_G)(1_{M(A_p)}\otimes d_G^s(\lambda (f))) 
\]
for all $a_p$ in $A_p$ and for all $f$ in $C_c(G)$. It is not difficult to
check that $\{\widehat{\alpha }^{(p)}\}_{p\in S(A)}$ is an inverse system of
injective $S$-morphisms and $\widehat{\alpha }=\lim\limits_{\stackunder{p\in
S(A)}{\leftarrow }}\widehat{\alpha }^{(p)}$ is a coaction of $G$ on $A\times
_{\alpha ,r}G$ which verifies the condition \label{*} $(*)$.

$(2).$ By Theorem 5.2.6 (ii) of [{\bf 13}], for each $\delta \in \Delta $
there is a continuous action $\widehat{\beta }^{(\delta )}$ of $G$ on $%
A_{\delta ^{}}\times _{\beta ^{(\delta )}}G$ such that 
\[
\widehat{\beta }^{(\delta )}(\beta ^{(\delta )}(a_{\delta
^{}})(1_{M(A_{\delta ^{}})}\otimes f))=(\beta ^{(\delta )}(a_{\delta
^{}})\otimes 1_G)(1_{M(A_{\delta ^{}})}\otimes (\overline{\text{id}%
_{C_0(G)}\otimes j_G^a})d_G^a(f)) 
\]
for all $a_{\delta ^{}}$ in $A_{\delta ^{}}$ and for all $f$ in $C_c(G)$.
Using this relation and Remark 5.3 it is not difficult to check that $\{%
\widehat{\beta }^{(\delta )}\}_{\delta \in \Delta }$ is an inverse system of
injective $S$-morphisms. Let $\widehat{\beta }=\lim\limits_{\stackunder{%
\delta \in \Delta }{\leftarrow }}\widehat{\beta }^{(\delta )}$. Then $%
\widehat{\beta }$ is a continuous action of $G$ on $A\times _{\beta _{}}G$
and moreover, it verifies the condition $(**)$.%
\endproof%

\smallskip

\smallskip The following theorem is a version of the Takai duality theorem
for inverse limit actions of a locally compact group on a locally $C^{*}$%
-algebra.

\begin{Theorem}
Let $G$\ be a locally compact group, let $A$\ be a locally $C^{*}$-algebra
and let $\alpha $\ be an inverse limit action of $G$\ on $A$. Then there is
an isomorphism $\Pi $\ from $A\otimes {\cal K}(L^2(G))$\ onto $\left(
A\times _{\alpha ,r}G\right) \times _{\widehat{\alpha }}G$\ such that 
\[
\widehat{\widehat{\alpha }}\circ \Pi =(\overline{\Pi \otimes \text{id}%
_{C_0(G)}})\circ (\alpha \otimes \text{ad}\rho ) 
\]
where $\rho $\ is the right regular representation of $L^1(G)$.
\end{Theorem}

\smallskip

\proof%
By Proposition 3.2 of [{\bf 10}], 
\[
A\otimes {\cal K}(L^2(G))=\lim\limits_{\stackunder{p\in S(A)}{\leftarrow }%
}A_p\otimes {\cal K}(L^2(G)) 
\]
up to an isomorphism of locally $C^{*}$-algebras

Since $\alpha $ is an inverse limit action, according to the proof of
Proposition 5.5 (1), 
\[
\widehat{\alpha }=\lim\limits_{\stackunder{p\in S(A)}{\leftarrow }}\widehat{%
\alpha }^{(p)} 
\]
where $\widehat{\alpha }^{(p)}$ is the dual coaction associated to $\alpha
^{(p)}$ for each $p\in S(A)$. Then, since the connecting maps of the inverse
system $\{A_p\times _{\alpha ^{(p)},r}G\}_{p\in S(A)}$ are all surjective,
by Proposition 5.5 (2), 
\[
\widehat{\widehat{\alpha }}=\lim\limits_{\stackunder{p\in S(A)}{\leftarrow }}%
\widehat{\widehat{\alpha }}^{(p)} 
\]
and by Remark 5.3, 
\[
\left( A\times _{\alpha ,r}G\right) \times _{\widehat{\alpha }%
}G=\lim\limits_{\stackunder{p\in S(A)}{\leftarrow }}\left( A_p\times
_{\alpha ^{(p)},r}G\right) \times _{\widehat{\alpha }^{(p)}}G 
\]
up to an isomorphism of locally $C^{*}$ -algebras.

Let $p\in S(A)$. According to Theorem 5.2 of [{\bf 13}], there is an
isomorphism $\Pi ^{(p)}$ from $A_p\otimes {\cal K}(L^2(G))$ onto $\left(
A_p\times _{\alpha ^{(p)},r}G\right) \times _{\widehat{\alpha }^{(p)}}G$
such that 
\[
\widehat{\widehat{\alpha }}^{(p)}\circ \Pi ^{(p)}=(\overline{\Pi
^{(p)}\otimes \text{id}_{C_0(G)}})\circ (\alpha ^{(p)}\otimes \text{ad}\rho
). 
\]
Moreover, 
\[
\Pi ^{(p)}(\alpha ^{(p)}(a_p)(1_{M(A_p)}\otimes \lambda (f)h))=\widehat{%
\alpha }^{(p)}(\alpha ^{(p)}(a_p)(1_{M(A_p)}\otimes \lambda
(f)))(1_{M(A_p)}\otimes 1_G\otimes h) 
\]
and 
\[
\Pi ^{(p)}((1_{M(A_p)}\otimes \lambda (f)h)\alpha ^{(p)}(a_p))=\widehat{%
\alpha }^{(p)}((1_{M(A_p)}\otimes \lambda (f))\alpha
^{(p)}(a_p))(1_{M(A_p)}\otimes 1_G\otimes h) 
\]
for all $f$ and $h$ in $C_c(G)$ and for all $a_p$ in $A_p$. Using these
relations and the fact that $A_p\otimes {\cal K}(L^2(G))$ is the $C^{*}$%
-subalgebra of $M(A_p\otimes {\cal K}(L^2(G)))$ generated by $\{\alpha
^{(p)}(a_p)(1_{M(A_p)}\otimes \lambda (f)h),$ $(1_{M(A_p)}\otimes \lambda
(f)h)\alpha ^{(p)}(a_p);f,h\in C_c(G),a_p\in A_p\}$, see Lemma 5.2.10 of [%
{\bf 13}], it is not difficult to check that $\{\Pi ^{(p)}\}_{p\in S(A)}$ is
an inverse system of $C^{*}$-isomorphisms.

Let $\Pi =\lim\limits_{\stackunder{p\in S(A)}{\leftarrow }}\Pi ^{(p)}$.
Then, clearly $\Pi $ is an isomorphism of locally $C^{*}$-algebras from $%
A\otimes {\cal K}(L^2(G))$ onto $\left( A\times _{\alpha ,r}G\right) \times
_{\widehat{\alpha }}G$ which satisfies the condition 
\[
\widehat{\widehat{\alpha }}\circ \Pi =(\overline{\Pi \otimes \text{id}_{C(G)}%
})\circ (\alpha \otimes \text{ad}\rho ) 
\]
and the theorem is proved.%
\endproof%

\smallskip

Since any action of a compact group on a locally $C^{*}$-algebra is an
inverse limit action, we have:

\begin{Corollary}
Let $G$ be a compact group, let $A$ be a locally $C^{*}$-algebra and let $%
\alpha $ be a continuous action of $G$ on $A$. Then there is an isomorphism $%
\Pi $\ from $A\otimes {\cal K}(L^2(G))$\ onto $\left( A\times _{\alpha
,r}G\right) \times _{\widehat{\alpha }}G$\ such that 
\[
\widehat{\widehat{\alpha }}\circ \Pi =(\overline{\Pi \otimes \text{id}%
_{C_0(G)}})\circ (\alpha \otimes \text{ad}\rho ) 
\]
where $\rho $\ is the right regular representation of $L^1(G)$.
\end{Corollary}

\smallskip

{\bf Acknowledgment}. The author is grateful to the referee for several
suggestion that improved the presentation of the paper.

\ {\sc Department of Mathematics, Faculty of Chemistry, University of
Bucharest, Bd. Regina Elisabeta nr.4-12, Bucharest, Romania\ \ }

{\bf e-mail address:\ \ mjoita@fmi.unibuc.ro}

\end{document}